\UseRawInputEncoding
\documentclass{article}

\usepackage{stmaryrd}
\usepackage{lmodern}
\usepackage{diagbox}
\usepackage{mathtools}
\usepackage{amssymb,amsmath}
\usepackage{graphicx}

\usepackage{mathrsfs}

\usepackage{color}
\usepackage{amsfonts}
\usepackage{epsfig}
\usepackage{graphics}

\usepackage{lscape}
\usepackage{latexsym}

\def\proof{\noindent {\bf Proof }}
\def\alp{{\alpha}}

\def\gam{{\gamma}}
\def\qed{~~\vrule height8pt width4pt depth0pt}

\def\Ub{{\bar U}}
\def\Vb{{\bar V}}
\def\Wb{{\bar W}}
\def\Xb{{\bar X}}
\def\Hb{{\bar H}}

\def\F2{{\mathbb F}_2}

\def\eps{{\varepsilon}}

\def\Scal{{\cal S}}

\def\Mcal{{\cal M}}
\def\Ecal{{\cal E}}
\def\RS{{\cal RS}}

\def\OT{\mathbb{O}_1}

\newtheorem{prop}{Proposition}

\newtheorem{lemma}{Lemma}
\newtheorem{thm}{Theorem}
\newtheorem{cor}{Corollary}

\newcommand{\fp}{{\mathbb F}_{p}}
\newcommand{\fq}{{\mathbb F}_{q}}
\newcommand{\Fq}{{\mathbb F}_{q}}

\newcommand{\va}{{\vec a}}

\newcommand{\Nt}{{ N}}
\newcommand{\qone}{n}

\begin{document}

\title {Counting polynomials over finite fields with \\
prescribed leading coefficients and linear factors}
 \author{
Zhicheng Gao\\
School of Mathematics and Statistics\\
Carleton University\\
Ottawa, Ontario\\
Canada K1S5B6\\
Email:~zgao@math.carleton.ca }
\maketitle

\begin{abstract}
We count the number of polynomials over finite fields with  prescribed leading coefficients and a given number of linear factors. This is equivalent to
counting codewords in Reed-Solomon codes which are at a certain distance from a received word.
We first apply the generating function approach, which is recently developed by the author and collaborators,
to derive expressions for the number of monic polynomials with prescribed leading coefficients and linear factors.
We then apply Li and Wan's sieve formula to simplify the expressions in some special cases.
Our results extend and improve some recent results  by Li and Wan, and Zhou, Wang and Wang.

\end{abstract}

\section{ Introduction}

Let $\fq$ be the finite field with $q$ elements where $q$ is a prime power. Let $D=\{x_1,x_2,\ldots, x_n\}$ be a given set of $n$ elements from $\fq$ where $n\le q$.
The Reed-Solomon code $\RS_{n,k}$ is the set of all
words of the form $(f(x_1),f(x_2),\ldots, f(x_n))$, where $f$ is a polynomial over $\fq$ of degree at most $k-1$.
 The code $\RS_{q,k}$ is called the standard Reed-Solomon code.
Recall the (Hamming) distance between two words $u,v$ is the number of non-zero entries in $u-v$. For the Reed-Solomon code,
the distance between two codewords represented by polynomials $f$ and $g$ is the number of distinct roots of $f-g$ in $D$, or equivalently the number of distinct
linear factors $x-\alp$ with $\alp\in D$.  Given a monic polynomial $f$ of degree $k+\ell$,
the number of codewords in $\RS_{n,k}$ whose distance from $f$ is exactly $n-r$ is equal to the
number of monic polynomials $g$ of degree $k+\ell$
such that $g$ has the same $\ell$ leading coefficients as $f$, and $g$ has exactly $r$ distinct roots in $D$.
Following the notation in \cite{LiWan20,ZhouWang17}, we shall denote this number by $N(f,r)$. Readers are referred to \cite{ChenMu07,LiWan20} for more information about $\RS_{n,k}$ and its connections to polynomials with prescribed leading coefficients and linear factors.

In \cite{LiWan20}, Li and Wan give a bound on $N(f,r)$ for the $\RS_{q,k}$ code using Weil's bound on character sums.  In  \cite{ZhouWang17}
simple explicit expressions are given for the $\RS_{q,k}$ code when $\ell=1$ and when $\ell=2$ with $q$ being an even power of an odd prime.
In this paper we extend  the results in \cite{ZhouWang17} to general $\RS_{n,k}$ codes.  We also show that our bound improves Li and Wan's bound  \cite[Theorem~1.5]{LiWan20}
for two prescribed coefficients.
Our approach follows that of \cite{GKW21b} using generating functions. The coefficients of the generating functions are from the group algebra generated by
equivalence classes consisting of polynomials with prescribed leading coefficients.

The remaining paper is organized as follows.
In Section~2 we set up relevant generating functions and use them to derive expressions for $N(f,r)$ for general $D\subseteq \fq$.
This generating function approach is developed recently in \cite{GKW21a, GKW21b}.
In Section~3 we derive explicit expressions for the case of prescribing one leading coefficient.
The result extends  \cite[Theorem~3.1]{ZhouWang17} to allow $D\cup \{0\}$ to be any additive subgroup of $\fq$.
In Section~4 we derive explicit expressions for the case of prescribing two leading coefficients. Exact expressions are obtained when $n$ is an even power of $p$, which extend  \cite[Theorem~4.3]{ZhouWang17} in two ways: first, $D$ can be  any subfield of $\fq$;
and second, the two leading coefficients are arbitrary (not just 0).
Asymptotic expressions are obtained for $N(f,r)$ when $n$ is any power of $p$.
We demonstrate that our error bound is exponentially smaller than the one given in \cite{LiWan20} when the degree of $f$ is equal to $p-1$ and $q=p^2$ for some prime $p$.
In Section~5 we give simple expressions for the expected value and variance of the distance between a given received word  and a random codeword in $\RS_{n,k}$. Section~6 concludes our paper.

\section{Generating functions for polynomials with prescribed leading coefficients and linear factors}

 Fix a positive integer $\ell$.
Given a polynomial  $f=x^d+f_1x^{d-1}+\cdots f_{d-1}x+f_d$, we shall call $f_1,\ldots, f_{\ell}$ the leading coefficients of $f$.
When we read the leading coefficients from left to right, missing coefficients are interpreted as zero.
Thus the leading coefficients of $f$
are the same as those of $x^kf$ for any $k\ge 0$.

 Let $\Mcal$ denote the set of monic polynomials  over $\Fq$,  $\Mcal_d$ be the subset of $\Mcal$ consisting of the monic polynomials of degree $d$, and let $\deg(f)$ denote the degree of a polynomial $f$.
Two polynomials $f,g\in \Mcal$ are said to be {\em equivalent} if they have the same $\ell$ leading coefficients. 

We shall use $\langle f\rangle_{\ell}$ to denote the equivalence class containing $f$,
where the subscript $\ell$ is usually omitted when it is clear from the context. For typographical convenience, we shall use $\langle a_1,\ldots, a_{\ell}\rangle$
to denote the equivalence class $\langle x^{\ell}+a_1x^{\ell-1}+\cdots +a_{\ell} \rangle$.
We also use $\Ecal_{\ell}$ (or simply $\Ecal$) to denote the set of all equivalence classes. Given $\eps\in \Ecal$, we shall also use $\Mcal_d(\eps)$ to denote the set of polynomials in $\Mcal_d$ which are  equivalent to $\eps$.

The following proposition shows that $\Ecal$, under the usual polynomial multiplication, forms a group. Its proof can be found  in \cite{EH91,Fom96, GKW21b, Hay65}.

\begin{prop}\label{prop:one}
For each given positive integer $\ell$, the set $\Ecal$ under the multiplication $\langle f\rangle \langle g\rangle =\langle fg\rangle$
forms a group  with $\langle 1\rangle$ being the identity element. Moreover,
\begin{itemize}
\item[(a)] Each  equivalence class is represented by a unique polynomial in $\Mcal_{\ell}$, and $|\Ecal|=q^{\ell}$.
\item[(b)] For each $d\ge \ell$ and each $\eps\in \Ecal$, there are exactly $q^{d-\ell}$ polynomials in $\Mcal_d(\eps)$.
\end{itemize}
\end{prop}

We shall use $0$ to denote the zero element of the group algebra $\mathbb{Q}\Ecal$ generated by the group $\Ecal$ over the field of rational numbers.

For a subset $S$ of $D$, let $\Nt_d(\eps, S)$ be the number of polynomials $f(x)$ in $\Mcal_d(\eps)$  such that $f(x)=g(x)\prod_{\alp\in S}(x+\alp)$ and $g(x)$
contains no linear factor $(x+\beta)$ for any $\beta\in \fq \setminus S$. Define $M_d(\eps, S)$ similarly by dropping the latter condition.
In other words, $N_d(\eps,S)$ counts polynomials with linear factors {\em exactly} in $S$, and  $M_d(\eps,S)$ counts polynomials with linear factors in $S$ and {\em possibly more}.

Given a polynomial $f\in \Mcal_{d}$, we follow the notation in \cite{LiWan20} and define
\begin{align*}
M(f,r)&=\sum_{S\subseteq D, |S|=r} M_{d}(\langle f\rangle,S),\\
N(f,r)&=\sum_{S\subseteq D, |S|=r} \Nt_{d}(\langle f\rangle,S).
\end{align*}
 Thus $N(f,r)$  counts polynomials in
$\Mcal_{d}(\langle f\rangle)$ with exactly $r$ linear factors associated with $D$.
We also note that $N(f,r)$ and $M(f,r)$ are related by the following sieve formula (see, e.g., \cite[Section~4.2]{Wilf}):
\begin{align}\label{SIEVE}
N(f,r)=\sum_{j=r}^{d}(-1)^{j-r}{j\choose r}M(f,j).
\end{align}
Hence we may either deal with the generating function of $N(f,r)$ directly or use the generating function of $M(f,r)$ and \eqref{SIEVE} to obtain the generating function of $N(f,r)$.

Now we define the following generating functions:
\begin{align}
F(z)&=\sum_{f\in \Mcal} \langle f\rangle z^{\deg(f)}= \sum_{d\ge 0}z^d\sum_{f\in \Mcal_d} \langle f\rangle ,\nonumber\\
F(z; S)&=\sum_{d\ge 0}z^d \sum_{\eps\in \Ecal} \Nt_d(\eps, S)\eps ,\label{eq:FN}\\
F(z; \supseteq\!\!S)&=\sum_{d\ge 0}z^d\sum_{\eps\in \Ecal} M_d(\eps, S)\eps . \nonumber
\end{align}
We note that $F(z)$, $F(z;S)$, and $F(z; \supseteq\!\! S)$ are formal power series with coefficients in the group algebra ${\mathbb Q}\Ecal$.

 As shown in \cite{GKW21b}, it is convenient to introduce the following two elements in ${\mathbb Q}\Ecal$.
\begin{align}
 E&=\frac{1}{q^{\ell}}\sum_{\eps\in \Ecal}\eps, \quad  J=\langle 1\rangle -E.\label{eq:E}
\end{align}

\begin{prop}\label{prop:GFun}
Let $E$ be defined in \eqref{eq:E}. Then
\begin{align}
F(z)&=\sum_{d=0}^{\ell-1}\sum_{f\in \Mcal_d}\langle f\rangle z^d+\frac{(qz)^{\ell}}{1-qz}E,\label{eq:F}\\
F(z;\supseteq\!\! S)&=\left(\prod_{\alp\in S}\langle x+\alp\rangle\right)z^{|S|} F(z),\label{eq:FSup}\\
F(z;S)&=\left(\prod_{\alp\in S}\langle x+\alp\rangle\right) z^{|S|}F(z)\left(\prod_{\beta\in \fq\setminus S}(\langle 1\rangle-\langle x+\beta\rangle z) \right). \label{eq:FS}
\end{align}
\end{prop}
\proof  Using Proposition~\ref{prop:one}, we have
\begin{align*}
F(z)&=\sum_{d=0}^{\ell-1}\sum_{f\in \Mcal_d}\langle f\rangle z^d+\sum_{d\ge \ell}\sum_{f\in \Mcal_d}\langle f\rangle z^d\\
&=\sum_{d=0}^{\ell-1}\sum_{f\in \Mcal_d}\langle f\rangle z^d+\sum_{d\ge \ell}\sum_{\eps\in \Ecal}q^{d-\ell}\eps z^d\\
&=\sum_{d=0}^{\ell-1}\sum_{f\in \Mcal_d}\langle f\rangle z^d+\frac{(qz)^{\ell}}{1-qz}E,
\end{align*}
which is \eqref{eq:F}.

Next we note that each monic polynomial containing the factor $\prod_{\alp\in S}(x+\alp)$ can be written uniquely as $f\prod_{\alp\in S}(x+\alp)$ for some monic polynomial $f$. So we have
\begin{align*}
F(z;\supseteq \!\!S)&=\sum_{d\ge 0}\sum_{f\in \Mcal_d}\left\langle f\prod_{\alp\in S}(x+\alp)\right\rangle z^{d+|S|}\\
&= z^{|S|}\prod_{\alp\in S}\langle x+\alp\rangle \sum_{d\ge 0}\sum_{f\in \Mcal_d}\langle f\rangle z^{d},
\end{align*}
which gives \eqref{eq:FSup}.

Finally we prove \eqref{eq:FS}. Let $\Mcal'$ denote the set of monic polynomials which does not contain any linear factor. Then each monic polynomial $f$ can be factored uniquely as
$f=g\prod_{\beta\in \fq}(x+\beta)^{d(\beta)}$ for some $g\in \Mcal'$ and non-negative integers $d(\beta)$. Thus
\begin{align*}
F(z)&=\left(\prod_{\beta\in \fq}\sum_{d(\beta)\ge 0}\langle x+\beta\rangle^{d(\beta)}z^{d(\beta)}\right)\sum_{g\in \Mcal'}\langle  g \rangle z^{\deg(g)}\\
&=\left(\prod_{\beta\in \fq}\frac{1}{\langle 1\rangle-z\langle x+\beta\rangle} \right)\sum_{g\in \Mcal'}\langle  g \rangle z^{\deg(g)},\\
F(z;S)&=\left(\prod_{\alp\in S}\sum_{d(\alp)\ge 1}\langle x+\alp\rangle^{d(\alp)}z^{d(\alp)}\right)\sum_{g\in \Mcal'}\langle  g \rangle z^{\deg(g)}\\
&=\left(\prod_{\alp\in S}\frac{z\langle x+\alp\rangle}{\langle 1\rangle-z\langle x+\alp\rangle} \right)\sum_{g\in \Mcal'}\langle  g \rangle z^{\deg(g)},
\end{align*}
and \eqref{eq:FS} follows immediately. \qed

The next proposition will be used to extract the coefficients from our generating functions.
\begin{prop}\label{prop:EJ} Let $E$ and $J$ be defined in \eqref{eq:E}. The following is true.
\begin{itemize}
\item[(a)] For each  $\eps\in \Ecal$, $E \eps =E$.
\item[(b)] $E^2=E$, $J^2=J$, and $EJ=0$.
\item[(c)] For each $d\ge \ell$, $\displaystyle J\sum_{f\in \Mcal_d}\langle f\rangle =0$.
\end{itemize}
\end{prop}
\proof (a) Noting that $\Ecal$ is a group, we have
\begin{align*}
E\eps&=\frac{1}{q^{\ell}}\sum_{\eta\in \Ecal}\eta\eps=\frac{1}{q^{\ell}}\sum_{\eta'\in \Ecal}\eta'=E.
\end{align*}
(b) It follows from part~(a) and Proposition~\ref{prop:one}(a) that
\begin{align*}
E^2&=\frac{1}{q^{\ell}}\sum_{\eta\in \Ecal}E\eta=\frac{1}{q^{\ell}}\sum_{\eta'\in \Ecal}E=E,\\
EJ&=E(\langle 1\rangle -E)=E-E^2=0,\\
J^2&=J(\langle 1\rangle -E)=J-JE=J.
\end{align*}
(c) Using part~(b) and Proposition~\ref{prop:one}(b), we obtain
\begin{align*}
J\sum_{f\in \Mcal_d}\langle f\rangle &=J\sum_{\eps\in \Ecal}q^{d-\ell}\eps=q^dJE=0.
\end{align*}
\qed

We now use the above propositions to obtain expressions for $N(f,r)$ and $M(f,r)$ by extracting coefficients from $F(z;S)$. In the rest of the paper, we shall use $[z^d\eps]F$ to denote the coefficient of $z^d\eps$ in the generating function $F$. It is  convenient to use the Iverson's bracket $\llbracket P\rrbracket$, which has value 1 if predicate $P$ is true and 0 otherwise.
It is also convenient to introduce the notations $X_m=D^m$ and $\Xb_m=\{(x_1,\ldots,x_m)\in X_m:x_i\ne x_j \hbox{ if }i\ne j \}$. That is, $\Xb_m$ is a subset of $X_m$ with distinct coordinates. The subscript $m$ is often omitted when it is clear from the context.
\begin{thm}\label{thm:zero}
\begin{align}
M(f,r)
&=\llbracket r\le k\rrbracket q^{k-r}{n\choose r}\nonumber\\
&~~+\frac{\llbracket k<r\le k+\ell\rrbracket}{r!}\sum_{ g\in \Mcal_{k+\ell-r},(x_1,\ldots,x_r)\in \Xb}\left\llbracket \langle g\rangle\prod_{j=1}^r\langle x+x_j\rangle=\langle f\rangle \right\rrbracket,\label{eq:M}\\
N(f,r)
&= q^{k-r}{n\choose r}\sum_{j=0}^{k-r}(-q)^{-j}{n-r\choose j}\nonumber\\
&~~+\sum_{i=1}^{\ell}(-1)^{k+i-r}{k+i\choose r}\frac{1}{(k+i)!}\sum_{ g\in \Mcal_{\ell-i},(x_1,\ldots,x_{k+i})\in \Xb}\left\llbracket \langle g\rangle\prod_{j=1}^{k+i}\langle x+x_j\rangle=\langle f\rangle \right\rrbracket.\label{eq:N}
\end{align}
\end{thm}
\proof Using Propositions~\ref{prop:GFun} and \ref{prop:EJ}, we obtain
\begin{align*}
EF(z;\supseteq \!\!S)&=E\left(\prod_{\alp\in S}\langle x+\alp\rangle\right) z^{|S|}\left(\sum_{d=0}^{\ell-1}\sum_{g\in \Mcal_d }Ez^d+E\frac{(qz)^{\ell}}{1-qz}\right)\\ 
&=z^{|S|}E\left(\sum_{d=0}^{\ell-1}q^dz^d+\frac{(qz)^{\ell}}{1-qz}\right),\\
JF(z;\supseteq \!\!S)&=z^{|S|}J\left(\prod_{\alp\in S}\langle x+\alp\rangle\right)\left(J\sum_{d=0}^{\ell-1}\sum_{g\in \Mcal_d}\langle g\rangle z^d\right) \nonumber\\
&=z^{|S|}(\langle 1\rangle-E)\left(\prod_{\alp\in S}\langle x+\alp\rangle\right)\left(\sum_{d=0}^{\ell-1}\sum_{g\in \Mcal_d}\langle g\rangle z^d\right) \nonumber\\
&=z^{|S|}\left(\prod_{\alp\in S}\langle x+\alp\rangle\right)\left(\sum_{d=0}^{\ell-1}\sum_{g\in \Mcal_d}\langle g\rangle z^d\right) -z^{|S|}E\left(\sum_{d=0}^{\ell-1}q^dz^d\right).
\end{align*}
Extracting the coefficients, we obtain
\begin{align*}
M_{k+\ell}(\eps,S)&=\left[z^{k+\ell}\eps\right]EF(z;\supseteq \!\!S)+ \left[z^{k+\ell}\eps\right]JF(z;\supseteq \!\!S)   \\
&=\left[z^{k+\ell}\eps\right]z^{|S|+\ell}q^{\ell}E\frac{1}{1-qz}+
\left[z^{k+\ell}\eps\right]z^{|S|}\left(\prod_{\alp\in S}\langle x+\alp\rangle\right)\left(\sum_{d=0}^{\ell-1}\sum_{g\in \Mcal_d}\langle g\rangle z^d\right)\nonumber\\
&=\left\llbracket k\ge |S| \right\rrbracket q^{k-|S|}+\left[z^{k+\ell-|S|}\eps\right]\left(\prod_{\alp\in S}\langle x+\alp\rangle\right)\left(\sum_{d=0}^{\ell-1}\sum_{g\in \Mcal_d}\langle g\rangle z^d\right)\nonumber\\
&=\llbracket |S|\le k\rrbracket q^{k-|S|}+\llbracket k<|S|\le k+\ell\rrbracket
\left \llbracket \exists g\in \Mcal_{k+\ell-|S|},\eps=\langle g\rangle\prod_{\alp\in S}\langle x+\alp\rangle \right \rrbracket.
\end{align*}
Setting $\eps=\langle f\rangle$, summing over all $S\subseteq D$ with $|S|=r$, and noting that each subset $\{x_1,\ldots,x_r\}$ corresponds to $r!$ ordered $r$-tuples in $\Xb$, we obtain \eqref{eq:M}. Substituting \eqref{eq:M} into \eqref{SIEVE}, we obtain
\begin{align}
N(f,r)&=\sum_{j=r}^{k+\ell}(-1)^{j-r}{j\choose r}\llbracket j\le k\rrbracket q^{k-j}{n\choose j} \nonumber\\
&~~+\sum_{j=r}^{k+\ell}(-1)^{j-r}{j\choose r}\frac{\llbracket k<j\le k+\ell\rrbracket}{j!}\sum_{ g\in \Mcal_{k+\ell-j},(x_1,\ldots,x_j)\in \Xb}\left\llbracket \langle g\rangle\prod_{i=1}^j\langle x+x_i\rangle=\langle f\rangle \right\rrbracket. \label{eq:Ntemp} 
\end{align}
Changing the summation index $j:=j+r$ and using the identity
\begin{align}\label{eq:BinIden}
{m\choose r}{n\choose m}={n\choose r}{n-r\choose m-r},
\end{align} 
we can rewrite the first line of \eqref{eq:Ntemp} as the first line of \eqref{eq:N}. The second line of \eqref{eq:Ntemp} is converted into the second line of \eqref{eq:N} by changing the summation index $j:=k+i$. \qed

Although it is less convenient to deal with $N_{k+\ell}(\eps,S)$ directly when $\ell\ge 2$, we note that the same argument can be used to derive an expression for $N_{k+\ell}(\eps,S)$ directly from $F(z;S)$, as shown below.
\begin{align*}
EF(z;S)&=\left(\prod_{\alp\in S}E\langle x+\alp\rangle\right) z^{|S|}EF(z)\left(\prod_{\beta\notin S}(E-E\langle x+\beta \rangle z) \right) \\
&= z^{|S|}E\left(\sum_{d=0}^{\ell-1}q^dz^d+\frac{(qz)^{\ell}}{1-qz}\right)(1-z)^{q-|S|},\\
JF(z;S)&=z^{|S|}J\left(\prod_{\alp\in S}\langle x+\alp\rangle\right)\left(\sum_{d=0}^{\ell-1}\sum_{f\in \Mcal_d}\langle f\rangle z^d\right)\prod_{\beta\notin S}(\langle 1\rangle-\langle x+\beta\rangle z) \\
&=z^{|S|}\left(\prod_{\alp\in S}\langle x+\alp\rangle\right)\left(\sum_{d=0}^{\ell-1}\sum_{f\in \Mcal_d}\langle f\rangle z^d\right)\prod_{\beta\notin S}(\langle 1\rangle-\langle x+\beta\rangle z)\\
&~~-z^{|S|}E\left(\sum_{d=0}^{\ell-1}q^dz^d\right)(1-z)^{q-|S|}.
\end{align*}

By \eqref{eq:FN}, we have
\begin{align}
N_{k+\ell}(\eps,S)&=\left[z^{k+\ell}\eps\right]F(z;S)\nonumber\\
&=\left[z^{k+\ell}\eps\right]\left(EF(z;S)+JF(z;S)\right)\nonumber\\
&=\left[z^{k+\ell}\eps\right]z^{|S|+\ell}q^{\ell}E\frac{1}{1-qz}(1-z)^{q-|S|}\nonumber\\
&~~~+\left[z^{k+\ell}\eps\right]z^{|S|}\left(\prod_{\alp\in S}\langle x+\alp\rangle\right)\left(\sum_{d=0}^{\ell-1}\sum_{g\in \Mcal_d}\langle g\rangle z^d\right)\prod_{\beta\in \fq\setminus S}(\langle 1\rangle-\langle x+ \beta\rangle z)\nonumber\\
&=\left[z^{k-|S|}\right]\frac{1}{1-qz}(1-z)^{q-|S|}\nonumber\\
&~~~+\left[z^{k+\ell-|S|}\eps\right]\left(\prod_{\alp\in S}\langle x+\alp\rangle\right)\left(\sum_{d=0}^{\ell-1}\sum_{g\in \Mcal_d}\langle g\rangle z^d\right)\prod_{\beta\in \fq\setminus S}(\langle 1\rangle-\langle x+\beta\rangle z)\nonumber\\
&=\sum_{j=0}^{k-|S|}{q-|S|\choose j}(-1)^jq^{k-|S|-j}\nonumber\\
&~~~+\left[z^{k+\ell-|S|}\eps\right]\left(\prod_{\alp\in S}\langle x+\alp\rangle\right)\left(\sum_{d=0}^{\ell-1}\sum_{g\in \Mcal_d}\langle g\rangle z^d\right)\prod_{\beta\in \fq\setminus S}(\langle 1\rangle-\langle x+\beta\rangle z).  \label{eq:NSexp}
\end{align}
When $\ell=0$, the second term in \eqref{eq:NSexp} disappears, and we have
\begin{align*}
\Nt_k(\eps,S)&=\sum_{j=0}^{k-|S|}{q-|S|\choose j}(-1)^jq^{k-|S|-j},\\
N\left(x^{k},r\right)&=\sum_{S\subseteq \fq,|S|=r}\Nt_k(\eps,S)= {q\choose r}q^{k-r}\sum_{j=0}^{k-r}{q-r\choose j}(-1)^jq^{-j}.
\end{align*}
This is \cite[Corollary 1.6]{LiWan20}, which is first obtained in \cite{KK90}.

\section{The case $\ell=1$}

In this section we treat the case $\ell=1$, that is,  monic polynomials with prescribed trace.

Setting $\ell=1$ in \eqref{eq:N} and noting $\Mcal_0=\{1\}$, we obtain
\begin{align*}
N(x^{k+1}+\gam x^{k},r)&=q^{k-r}{n\choose r}\sum_{j=0}^{k-r}(-q)^{-j}{n-r\choose j}\nonumber\\
&~~(-1)^{k+1-r}{k+1\choose r}\frac{1}{(k+1)!}\sum_{ (x_1,\ldots,x_{k+1})\in \Xb}\left\llbracket \prod_{j=1}^{k+1}\langle x+x_j\rangle=\langle x+\gam\rangle \right\rrbracket.
\end{align*}
Since $\ell=1$, we have
\[
\prod_{j=1}^{k+1}\langle x+x_j\rangle=\left \langle x+\sum_{j=1}^{k+1} x_j \right\rangle.
\]
It follows that
\begin{align}\label{eq:N1}
N(x^{k+1}+\gam x^{k},r)&=q^{k-r}{n\choose r}\sum_{j=0}^{k-r}(-q)^{-j}{n-r\choose j}
+(-1)^{k+1-r}{k+1\choose r}\frac{\Ub_{k+1}(\gam)}{(k+1)!},
\end{align}
where
\begin{align}\label{eq:Ub}
\Ub_m(\gam)=\sum_{ (x_1,\ldots,x_m)\in \Xb}\left\llbracket \sum_{j=1}^{m}x_j=\gam \right\rrbracket.
\end{align}

We now evaluate $\Ub_m(\gam)$ using  Li-Wan's sieve formula \cite[Theoem~3.1]{LiWan20}.  
Let $\Scal_m$ be the symmetric group on $\{1,2,\ldots,m\}$. For each $\tau\in \Scal_m$, define $X_{\tau}\subseteq X$ such that $x_i=x_j$ when $i,j$ belong to the same cycle of $\tau$.  Let $h$ be a complex-valued function defined on $X$, and define
\begin{align*}
\Hb=\sum_{(x_1,\ldots,x_m)\in \Xb}h(x_1,\ldots,x_m),~~H(\tau)=\sum_{(x_1,\ldots,x_m)\in {X_\tau}}h(x_1,\ldots,x_m).
\end{align*}
Theorem~3.1 of \cite{LiWan20} states
\begin{align}\label{eq:sieve}
\Hb=\sum_{\tau\in \Scal_m}(-1)^{m-l(\tau)}H(\tau),
\end{align}
where $l(\tau)$ denotes the number of cycles of $\tau$.
We shall need to make some assumptions on the set $D$, which is associated with the Reed-Solomon code $\RS_{n,k}$. In the rest of this section, we shall assume that either $D$ or $D\cup \{0\}$ is an additive subgroup of $\fq$.  Define
\begin{align} \label{eq:Line}
U_m(\gam)=\sum_{(x_1,\ldots,x_m)\in X}\left\llbracket\sum_{j=1}^m x_j=\gam \right\rrbracket.
\end{align}
It is clear that $U_r(\gam)=0$ when $\gam\notin D\cup \{0\}$. So we shall assume $\gam\in D\cup \{0\}$ below.

It is known (see \cite{GZ21} for general groups and \cite{GMW18,LiWan12} for abelian groups) that
\begin{align} \label{eq:Ur}
U_m(\gam)=\left\{\begin{array}{ll}
n^{m-1} & \hbox{ if $D$ is an additive subgroup of $\fq$},\\
\frac{n^m}{n+1}+\frac{(n+1)\llbracket \gam=0\rrbracket -1}{n+1}(-1)^m & \hbox{ if $0\notin D$ and $D\cup\{0\}$ is an additive group}.
\end{array}
\right.
\end{align}

Setting $h(x_1,\ldots,x_m)=\llbracket x_1+\cdots+x_m=\gam \rrbracket$, we may rewrite \eqref{eq:Ub} and \eqref{eq:Line} as
\begin{align*}
\Ub_m(\gam)&=\sum_{(x_1,\ldots,x_m)\in \Xb} h(x_1,\ldots,x_m),\\
U_m(\gam)&=\sum_{(x_1,\ldots,x_m)\in X} h(x_1,\ldots,x_m).
\end{align*}
Now \eqref{eq:sieve} becomes
\begin{align}\label{eq:Ub0}
\Ub_m(\gam)=\sum_{\tau\in \Scal_m}(-1)^{m-l(\tau)}U(\tau;\gam),
\end{align}
where
\begin{align}\label{eq:Utau0}
U(\tau;\gam)=\sum_{(x_1,\ldots,x_m)\in X_{\tau}} \llbracket x_1+\cdots+x_m=\gam \rrbracket.
\end{align}

In \eqref{eq:Utau0}, the sum of  $x_j$ over all $j$ in a cycle of length $t$ is equal to $tx_j$. Recall that $p$ is the characteristic of $\fq$. We note that $tx_j=0$ when $p\mid t$,  hence each cycle whose length is a multiple of $p$ contributes a factor $n$. When $p\nmid t$, $t$ has an inverse in $\fp$ and hence $tx_j$ can be treated as an independent variable.
Let $l'(\tau)$ be the number of cycles of $\tau$ whose length is not a multiple of $p$. We then have
\begin{align}\label{eq:Utau}
U(\tau;\gam)&=\llbracket l'(\tau)>0 \rrbracket \qone^{l(\tau)-l'(\tau)}U_{l'(\tau)}(\gam)+\llbracket l'(\tau)=0, \gam=0 \rrbracket  n^{l(\tau)}.
\end{align}

We are ready to prove the following result, which generalizes \cite[Theorem~3.1]{ZhouWang17}, allowing $D$ to be any additive subgroup of $\fq$.  We note that this result was also proved in \cite{LiWan12}. We include its proof here because the same argument is also used to prove the subsequent Theorems.
\begin{thm}\label{thm:one} Suppose $D$ is an additive subgroup of $\fq$. Then
\begin{align*}
N(x^{k+1}+\gam x^k,r)&=q^{k-r}{n\choose r}\sum_{j=0}^{k-r}(-q)^{-j}{n-r\choose j}
+\frac{\llbracket \gam\in D\rrbracket}{n}(-1)^{k+1-r}{n\choose r}{n-r\choose k+1-r}\nonumber\\
&~~+\frac{\llbracket p\mid (k+1),\gam\in D\rrbracket}{n}\left(n\llbracket \gam=0\rrbracket-1\right)(-1)^{r+(k+1)/p}{k+1\choose r}{n/p\choose (k+1)/p}.
\end{align*}
\end{thm}
\proof
Substituting the first expression of \eqref{eq:Ur} into \eqref{eq:Utau}, we obtain
\begin{align*}
U(\tau;\gam)
&=\llbracket l'(\tau)>0 \rrbracket \qone^{l(\tau)-1}+\llbracket l'(\tau)=0, \gam=0 \rrbracket  \qone^{l(\tau)}.
\end{align*}
It follows from \eqref{eq:Ub0} that
\begin{align}
\Ub_m(\gam)&=\sum_{\tau\in \Scal_m}(-1)^{m-l(\tau)}U(\tau;\gam)  \nonumber\\
&=\sum_{\tau\in \Scal_m}(-1)^{m-l(\tau)}\qone^{l(\tau)-1}+ \sum_{l'(\tau)=0,\tau\in \Scal_m} (-1)^{m-l(\tau)}\left(\llbracket \gam=0 \rrbracket \qone^{l(\tau)}- \qone^{l(\tau)-1}\right). \label{eq:Ubr}
\end{align}

To simplify the sums above, we recall the exponential generating function of permutations with respect to the total size and the cycle length.  Here we need to keep a separate count for those cycles whose lengths are not multiples of $p$.
The standard set construction (see, e.g., \cite{FlaSed09}) gives the following:
\begin{align}\label{eq:PGF}
\sum_{m\ge 0}\frac{z^m}{m!}\sum_{\tau\in \Scal_m}u^{l(\tau)}w^{l'(\tau)}
&=\exp\left(u\sum_{j\ge 1, p\mid j}z^j/j+uw\sum_{j\ge 1, p\nmid j}z^j/j\right)\nonumber\\
&=\exp\left((u-uw)\sum_{j\ge 1, p\mid j}z^j/j+uw\sum_{j\ge 1}z^j/j\right)\nonumber\\
&=\exp\left(\frac{u-uw}{p}\ln \frac{1}{1-z^p}+uw\ln \frac{1}{1-z}\right)\nonumber\\
&=(1-z)^{-uw}(1-z^p)^{-(u-uw)/p}.
\end{align}
At this stage, it is convenient to define
\begin{align}
A_m(u,w)&=[z^m](1-z)^{-uw}(1-z^p)^{-(u-uw)/p}\nonumber\\
&=\sum_{0\le j\le m/p}{uw+m-jp-1\choose m-jp}{(u-uw)/p+j-1\choose j}. \label{eq:Auw}
\end{align}
Thus \eqref{eq:PGF} can be rewritten as
\begin{align}
\frac{1}{m!}\sum_{\tau\in \Scal_m}u^{l(\tau)}w^{l'(\tau)}=A_m(u,w). \label{eq:Stirling3}
\end{align}

Setting $w=1$ and $w=0$ in \eqref{eq:Auw}, respectively, we obtain
\begin{align}
A_m(u,1)&=\frac{1}{m!}\sum_{\tau\in \Scal_m}u^{l(\tau)}=(-1)^m{-u\choose m}, \label{eq:Stirling1}\\
A_m(u,0)&=\frac{1}{m!}\sum_{\tau\in \Scal_m, l'(\tau)=0}u^{l(\tau)}=(-1)^{m/p}{-u/p\choose m/p}\llbracket p\mid m \rrbracket . \label{eq:Stirlingp}
\end{align}

Using \eqref{eq:Ubr}, \eqref{eq:Stirling1} and \eqref{eq:Stirlingp}, we obtain
\begin{align}
\frac{\Ub_m(\gam)}{m!}&= \frac{1}{n}(-1)^mA_m(-n,1)+\frac{n\llbracket \gam=0 \rrbracket -1}{n}(-1)^mA_m(-n,0)\nonumber\\
&=\frac{1}{\qone}{\qone\choose m}+(-1)^{m+m/p}\frac{n\llbracket \gam=0\rrbracket-1}{n}\llbracket p\mid m\rrbracket {n/p\choose m/p},
\end{align}
which is an extension of \cite[Lemma~3.1]{ZhouWang17}.  Substituting this into \eqref{eq:N1}, we obtain
\begin{align*}
 N(x^{k+1}+\gam x^{k},r)&=\sum_{j=r}^k(-1)^{j-r}{j\choose r} {n\choose j}q^{k-j}
+\frac{\llbracket\gamma\in D\rrbracket}{n}(-1)^{k+1-r}{k+1\choose r}{n\choose k+1} \\
&~+\llbracket \gamma\in D,  p\mid (k+1)\rrbracket  (-1)^{r+(k+1)/p}\frac{n\llbracket \gam=0\rrbracket-1}{n} {k+1\choose r}{n/p\choose (k+1)/p}.
\end{align*}
Changing the summation index  $j:=j-r$, and using \eqref{eq:BinIden},
we complete the proof. \qed

\medskip

When $0\notin D$, we have the following result.
\begin{thm}\label{thm:two} Suppose $0\notin D$ and $D\cup\{0\}$ is an additive subgroup of $\fq$.
\begin{itemize}
\item[(a)] If $\gam\in D$, then
\begin{align*}
N(x^{k+1}+\gam x^k,r)&=q^{k-r}{n\choose r}\sum_{j=0}^{k-r}(-q)^{-j}{n-r\choose j}
+\frac{1}{n+1}(-1)^{k+1-r}{n\choose r}{n-r\choose k+1-r}\nonumber\\
&~~-\frac{1}{n+1}(-1)^{r}{k+1\choose r}\sum_{0\le j\le (k+1)/p}(-1)^j{(n+1)/p\choose j}.
\end{align*}
\item[(b)] If $\gam\notin D\cup\{0\}$, then
\begin{align*}
N(x^{k+1}+\gam x^k,r)&=q^{k-r}{n\choose r}\sum_{j=0}^{k-r}(-q)^{-j}{n-r\choose j}.
\end{align*}
\item[(c)] If $\gam=0$, then
\begin{align*}
N(x^{k+1}+\gam x^k,r)&=q^{k-r}{n\choose r}\sum_{j=0}^{k-r}(-q)^{-j}{n-r\choose j}
+\frac{1}{n}(-1)^{k+1-r}{n\choose r}{n-r\choose k+1-r}\nonumber\\
&~~+\frac{n}{n+1}(-1)^{r}{k+1\choose r}\sum_{0\le j\le (k+1)/p}(-1)^j{(n+1)/p\choose j}.
\end{align*}
\end{itemize}
\end{thm}
\proof  The proof is similar to that of Theorem~\ref{thm:one}. We simply use the second expression of \eqref{eq:Ur} for $U_m(\gam)$. Now \eqref{eq:Utau} becomes
\begin{align*}
U(\tau;\gam)
&=\llbracket l'(\tau)>0 \rrbracket n^{l(\tau)-l'(\tau)}\left(\frac{n^{l'(\tau)}}{n+1}+\frac{(n+1)\llbracket \gam=0\rrbracket -1}{n+1}(-1)^{l'(\tau)}\right)\nonumber\\
&~~+\llbracket l'(\tau)=0, \gam=0 \rrbracket  n^{l(\tau)}.
\end{align*}
Consequently \eqref{eq:Ubr} becomes (note the cancelation  of the terms involving $\llbracket l'(\tau)=0\rrbracket$)
\begin{align*}
\frac{\Ub_m(\gam)}{m!}
&=\frac{(-1)^m}{n+1}(-1)^m\frac{1}{m!}\sum_{\tau\in \Scal_m}(-n)^{l(\tau)}+\frac{(-1)^m((n+1)\llbracket \gam=0\rrbracket -1)}{n+1} \frac{1}{m!}\sum_{\tau\in \Scal_m} (-n)^{l(\tau)-l'(\tau)}.
\end{align*}
Using \eqref{eq:Auw} again, we obtain
\begin{align*}
\frac{\Ub_m(\gam)}{m!}
&=\frac{1}{n+1}(-1)^mA_m(-n,1)+ \frac{(n+1)\llbracket \gam=0\rrbracket -1}{n+1}(-1)^mA_m(-n,-1/n)  \\
&=\frac{1}{n+1}{n\choose m}+\frac{(n+1)\llbracket \gam=0\rrbracket -1}{n+1}(-1)^m\sum_{0\le j\le m/p}(-1)^j{(n+1)/p\choose j}.
\end{align*}
Substituting this into \eqref{eq:N1}, we complete the proof. \qed

\section{The case $\ell=2$}
\def\fqone{{\mathbb F}_{n}}
\def\qone{{n}}

When $\ell\ge 2$, the expressions for $N(f,r)$ and $M(f,r)$ become much more complicated because they involve system of polynomial equations arising from Iverson's bracket in \eqref{eq:N}.  Throughout this section, we shall assume $\ell=2$ and $D=\fqone$ for some $n\mid q$.

Writing $f=x^{k+2}+\gam_1 x^{k+1}+\gam_2 x^{k}$, we obtain from \eqref{eq:N} that
\begin{align}
N(f,r)
&=q^{k-r}{n\choose r}\sum_{j=0}^{k-r}(-q)^{-j}{n-r\choose j}\nonumber\\
&~~+(-1)^{k+1-r}{k+1\choose r}\frac{1}{(k+1)!}\sum_{ \beta\in \fq,(x_1,\ldots,x_{k+1})\in \Xb}\left\llbracket \langle x+\beta\rangle \prod_{j=1}^{k+1}\langle x+x_j\rangle=\langle \gam_1,\gam_2\rangle \right\rrbracket\nonumber\\
&~~+(-1)^{k+2-r}{k+2\choose r}\frac{1}{(k+2)!}\sum_{ (x_1,\ldots,x_{k+2})\in \Xb}\left\llbracket \prod_{j=1}^{k+2}\langle x+x_j\rangle=\langle \gam_1,\gam_2\rangle \right\rrbracket\nonumber.
\end{align}
Since $\ell=2$, we have
\[
\prod_{j=1}^{m}\langle x+x_j\rangle=\left\langle \sum_{j=1}^m x_j,\sum_{1\le i<j\le m}x_ix_j \right\rangle.
\]
It follows that
\begin{align}
N(f,r)&=q^{k-r}{n\choose r}\sum_{j=0}^{k-r}(-q)^{-j}{n-r\choose j}\nonumber\\
&~~~+(-1)^{k+2-r}{k+2\choose r}\frac{\Vb_{k+2}(\gam_1,\gam_2) }{(k+2)!} \label{eq:N2}\\
&~~~+ (-1)^{k+1-r}{k+1\choose r}\frac{ \Wb_{k+1}(\gam_1,\gam_2)}{(k+1)!}, \nonumber
\end{align}
where
\begin{align}
\Vb_m(\gam_1,\gam_2)&=\sum_{(x_1,\ldots,x_m)\in \Xb}\left\llbracket\sum_{j=1}^m x_j=\gam_1,~~ \sum_{1\le i<j\le m} x_ix_j=\gam_2 \right\rrbracket, \label{eq:Vb} \\
\Wb_m(\gam_1,\gam_2)&=\sum_{y\in \fq,(x_1,\ldots,x_m)\in \Xb}\left\llbracket y+\sum_{j=1}^m x_j=\gam_1,~~ y\sum_{j=1}^m x_j+\sum_{1\le i<j\le m} x_ix_j=\gam_2 \right\rrbracket. \label{eq:Wb}
\end{align}

In the following, we shall assume $p>2$. Thus \eqref{eq:Vb} and \eqref{eq:Wb} are,  respectively, equivalent to the following diagonal systems:
\begin{align*}
\Vb_m(\gam_1,\gam_2)&=\sum_{(x_1,\ldots,x_m)\in \Xb}\left\llbracket\sum_{j=1}^m x_j=\gam_1,~~ \sum_{j=1}^m x^2_j=\gam_1^2-2\gam_2\right\rrbracket, \\
\Wb_m(\gam_1,\gam_2)&=\sum_{y\in \fq,(x_1,\ldots,x_m)\in \Xb}\left\llbracket y+\sum_{j=1}^m x_j=\gam_1,~~ y^2+\sum_{j=1}^m x^2_j=\gam_1^2-2\gam_2 \right\rrbracket.
\end{align*}

Define
\begin{align}\label{eq:V}
V_m(\gam_1,\gam_2)&=\sum_{(x_1,\ldots,x_m)\in X}\left\llbracket\sum_{j=1}^m x_j=\gam_1,~~ \sum_{j=1}^m x^2_j=\gam_1^2-2\gam_2\right\rrbracket.
\end{align}

To apply the sieve formula \eqref{eq:sieve}, we shall need to deal with the number of solutions to a slightly more general system:
\begin{align*}
V_m(\va;a_0,b_0)=\sum_{(x_1,\ldots,x_m)\in X}\left\llbracket \sum_{j=1}^m a_jx^2_j=a_0,~~\sum_{j=1}^m a_jx_j=b_0\right\rrbracket,
\end{align*}
where $a_1,\ldots, a_m\in \fp^*$.

When $q$  is an even power of an odd prime $p$,  the quadratic character $\eta$ over $\fq$ takes value 1 over $\fp^*$. Thus we have the following result, which is a special case of \cite[Lemma~4.1]{ZhouWang17}, written in a more compact form. \textcolor{red}{The quadratic character $\eta$ was missing in the published version \cite{pub22}.  This affacts Lemmas~1 and 2 and Theorem~4 in subsection 4.1, which are  indicated in red color}.

\begin{prop} \label{prop:three}
Assume  $q$ is an even power of an odd prime $p$, $D=\fqone$, $a_0,b_0\in D$, and let 
\[R_m(\va;a_0,b_0)=V_m(\va;a_0,b_0)-\qone^{m-2}.
\]
 Then
\[
\left|R_m(\va;a_0,b_0)\right|\le  \qone^{m/2}.
\]
Moreover, if $n$ is an even power of $p$, writing
$a=a_1\cdots a_m$  and  $b=a_1+\cdots+a_m$,  then we have
\begin{align*}
R_m(\va;a_0,b_0)&=\llbracket 2\mid m\rrbracket \qone^{(m-2)/2} \left(
\llbracket b_0=0, p\mid b \rrbracket (\qone\llbracket a_0=0\rrbracket-1)+\llbracket p\nmid b\rrbracket \textcolor{red}{\eta(b_0^2-ba_0)}\right)\\
&+\llbracket 2\nmid m\rrbracket \qone^{(m-1)/2}\left(
\llbracket  b_0=0,p\mid b\rrbracket \textcolor{red}{\eta(a_0)}+\llbracket p\nmid b \rrbracket 
\frac{\qone\llbracket b_0^2= ba_0\rrbracket-1}{\qone}\right).
\end{align*}
\end{prop}

\subsection{Exact expressions for the case  when $n$ is an even power of $p$}
In this subsection, we derive exact expressions for $N(x^{k+2}+\gamma_1x^{k+1}+\gamma_2x^k,r)$ when $n$ is an even power of $p$.
We first apply Proposition~\ref{prop:three} to derive the following result.
\begin{lemma} \label{lemma1}
Suppose $\qone$ is an even power of $p$, $D=\fqone$, and $\gam_1,\gam_2\in D$. Let $A_m(u,w)$ be defined by \eqref{eq:Auw}. Then the following holds.
\begin{itemize}
\item[(a)] For $p\mid m$,
\begin{align*}
\frac{\Vb_m(\gam_1,\gam_2)}{m!}
&=\frac{1}{\qone^2}\left({\qone\choose m}-{n/p\choose m/p}\right)+\frac{\llbracket \gam_1=0\rrbracket }{n} {n/p\choose m/p}\\
&~+\frac{\llbracket  \gam_1=0\rrbracket}{2n} \left(\qone\llbracket \gam_2=0\rrbracket-1+\textcolor{red}{\eta(\gamma_2)}\sqrt{n}   \right)(-1)^mA_m(-n,1/\sqrt{n})\\
&~+\frac{\llbracket  \gam_1=0\rrbracket}{2n} \left(\qone\llbracket \gam_2=0\rrbracket-1-\textcolor{red}{\eta(\gamma_2)}\sqrt{n}   \right)(-1)^mA_m(-n,-1/\sqrt{n}).
\end{align*}
\item[(b)] For $p\nmid m$,
\begin{align*}
\frac{\Vb_m(\gam_1,\gam_2)}{m!}
&=\frac{1}{\qone^2}{\qone\choose m}\\
&~+\frac{1}{2\qone\sqrt{n}}\left(\textcolor{red}{\eta((m-1)\gamma_1^2-2m\gamma_2)}\sqrt{n} +n\llbracket (m-1)\gam_1^2= 2m\gam_2\rrbracket -1\ \right)(-1)^m A_m(-n,1/\sqrt{n})\\
&~+\frac{1}{2\qone\sqrt{n}}\left(\textcolor{red}{\eta((m-1)\gamma_1^2-2m\gamma_2)}\sqrt{n} -n\llbracket (m-1)\gam_1^2= 2m\gam_2\rrbracket +1\ \right)(-1)^m A_m(-n,-1/\sqrt{n}).
\end{align*}
\end{itemize}
\end{lemma}

\proof  Define $V(\tau;\gam_1,\gam_2)$ similarly for those solutions to \eqref{eq:V} restricted to $X_{\tau}$.  Recall that $l'(\tau)$ is the number of cycles of $\tau$ whose lengths are not multiples of $p$. Let $\va$ be the vector of cycle lengths of $\tau$ which are not multiples of $p$. We have
\begin{align*}
V(\tau;\gam_1,\gam_2)&=\llbracket l'(\tau)>0\rrbracket n^{l(\tau)-l'(\tau)}V_{l'(\tau)}(\va;\gam_1^2-2\gam_2,\gam_1)
+\llbracket l'(\tau)=0,\gam_1=\gam_2=0\rrbracket n^{l(\tau)}.
\end{align*}

Noting
\begin{align*}
b=\sum_{j=1}^{l'(\tau)}a_j=m \pmod p,
\end{align*}
we see that conditions $p\mid b$ and $b_0^2=ba_0$ are equivalent to $p\mid m$ and $(m-1)\gam_1^2=2m\gam_2$, respectively. Thus we have
\begin{align}
\frac{\Vb_m(\gam_1,\gam_2)}{m!}&=\frac{1}{m!}\sum_{\tau\in \Scal_m,l'(\tau)>0}(-1)^{m-l(\tau)} n^{l(\tau)-l'(\tau)}V_{l'(\tau)}(\va;\gam_1^2-2\gam_2,\gam_1)\nonumber\\
&~~~+\frac{\llbracket \gam_1=\gam_2=0\rrbracket}{m!}\sum_{\tau\in \Scal_m,l'(\tau)=0}(-1)^{m-l(\tau)} n^{l(\tau)}\nonumber \\
&=\frac{1}{n^2}(-1)^m\left(A_m(-n,1)-A_m(-n,0)\right)
+\llbracket \gam_1=\gam_2=0\rrbracket (-1)^mA_m(-n,0)\nonumber\\
&~~~+\frac{1}{m!}\sum_{\tau\in \Scal_m,2\mid l'(\tau),l'(\tau)>0}(-1)^{m-l(\tau)} n^{l(\tau)-l'(\tau)}R_{l'(\tau)}(\va;\gam_1^2-2\gam_2,\gam_1)  \label{eq:Vbm}\\
&=\frac{1}{n^2}{n\choose m}+\left(\llbracket \gam_1=\gam_2=0\rrbracket -\frac{1}{n^2}\right)\llbracket p\mid m\rrbracket {n/p\choose m/p} \nonumber \\
&~~~+\frac{1}{m!}\sum_{\tau\in \Scal_m,2\mid l'(\tau),l'(\tau)>0}(-1)^{m-l(\tau)} n^{l(\tau)-l'(\tau)}R_{l'(\tau)}(\va;\gam_1^2-2\gam_2,\gam_1) \label{eq:Vb9} \\
&~~~+\frac{1}{m!}\sum_{\tau\in \Scal_m,2\nmid l'(\tau)}(-1)^{m-l(\tau)} n^{l(\tau)-l'(\tau)}R_{l'(\tau)}(\va;\gam_1^2-2\gam_2,\gam_1).  \nonumber
\end{align}

For part~(a), we use Proposition~\ref{prop:three} and \eqref{eq:Vb9} to obtain
\begin{align*}
\frac{\Vb_m(\gam_1,\gam_2)}{m!}&=\frac{1}{n^2}{n\choose m}+\left(\llbracket \gam_1=\gam_2=0\rrbracket -\frac{1}{n^2} \right){n/p\choose m/p}\nonumber\\
&~~+
\llbracket  \gam_1=0\rrbracket \frac{\qone\llbracket \gam_2=0\rrbracket-1}{\qone} (-1)^m
 \frac{1}{m!}\sum_{\tau\in \Scal_m,2\mid l'(\tau),l'(\tau)>0}(-n)^{l(\tau)} n^{-l'(\tau)/2}\nonumber\\
&~~+\eta(\gamma_2)\llbracket \gam_1=0\rrbracket   \frac{\sqrt{n}}{\qone}(-1)^m
 \frac{1}{m!}\sum_{\tau\in \Scal_m,2\nmid l'(\tau)}(-n)^{l(\tau)} n^{-l'(\tau)/2}.
\end{align*}
Using \eqref{eq:Stirling3} and \eqref{eq:Stirlingp}, and
\begin{align*}
\frac{1}{m!}\sum_{\tau\in \Scal_m,2\mid l'(\tau),l'(\tau)>0}(-n)^{l(\tau)} n^{-l'(\tau)/2}&=\frac{1}{2}(A_m(-n,1/\sqrt{n})+A_m(-n,-1/\sqrt{n})-A_m(-n,0)),\\
\frac{1}{m!}\sum_{\tau\in \Scal_m,2\nmid l'(\tau)}(-n)^{l(\tau)} n^{-l'(\tau)/2}&=\frac{1}{2}(A_m(-n,1/\sqrt{n})-A_m(-n,-1/\sqrt{n})),
\end{align*}
 we obtain
\begin{align*}
\frac{\Vb_m(\gam_1,\gam_2)}{m!}
&=\frac{1}{\qone^2}{\qone\choose m}+\frac{n\llbracket \gam_1=0\rrbracket -1}{n^2} {n/p\choose m/p}\\
&~+\llbracket  \gam_1=0\rrbracket  (-1)^m\frac{\qone\llbracket \gam_2=0\rrbracket-1}{2\qone} \left(A_m(-n,1/\sqrt{n})+A_m(-n,-1/\sqrt{n})\right)\\
&~~+\llbracket \gam_1=0\rrbracket  (-1)^m\textcolor{red}{\eta(\gamma_2)}\frac{\sqrt{n}}{2\qone} \left(A_m(-n,1/\sqrt{n})-A_m(-n,-1/\sqrt{n})\right).
\end{align*}
This completes the proof of part~(a).

\medskip

For part~(b), we use Proposition~\ref{prop:three}  and the same argument as above to obtain
\begin{align*}
\frac{\Vb_m(\gam_1,\gam_2)}{m!}
&=\frac{1}{\qone^2}{\qone\choose m}\\
&~+\textcolor{red}{\eta((m-1)\gamma_1^2-2m\gamma_2)}\frac{1}{2n} (-1)^m \left(A_m(-n,1/\sqrt{n})+A_m(-n,-1/\sqrt{n})\right)\\
&~~+\frac{\sqrt{n}}{2\qone^2}  \left(n\llbracket (m-1)\gam_1^2= 2m\gam_2\rrbracket -1\right) (-1)^m\left(A_m(-n,1/\sqrt{n})-A_m(-n,-1/\sqrt{n})\right).
\end{align*}
\qed

\medskip

A similar proof gives the following.
\begin{lemma} \label{lemma2}
Suppose $\qone$ is an even power of $p$, $D=\fqone$, and $\gam_1,\gam_2\in D$. Then the following holds.
\begin{itemize}
\item[(a)] For $p\mid (m+1)$,
\begin{align*}
\frac{\Wb_m(\gam_1,\gam_2)}{m!}
&=\frac{1}{n}{n\choose m}\\
&~+\frac{\llbracket \gam_1=0\rrbracket}{2\sqrt{n}}\left(\textcolor{red}{\eta(\gam_2)}\sqrt{n}+n\llbracket \gam_2=0\rrbracket-1\right)(-1)^mA_m(-n,1/\sqrt{\qone})\\
&~+\frac{\llbracket \gam_1=0\rrbracket}{2\sqrt{n}}\left(\textcolor{red}{\eta(\gam_2)}\sqrt{n}-n\llbracket \gam_2=0\rrbracket+1\right)(-1)^mA_m(-n,-1/\sqrt{\qone}).
\end{align*}
\item[(b)] For $p\nmid (m+1)$,
\begin{align*}
\frac{\Wb_m(\gam_1,\gam_2)}{m!}
&=\frac{1}{n}{n\choose m}\\
&~+\frac{1}{2n}\left(n\llbracket m\gam_1^2=2(m+1)\gam_2\rrbracket -1+\textcolor{red}{\eta(m\gam_1^2-2(m+1)\gam_2)}\sqrt{n}  \right)(-1)^mA_m(-n,1/\sqrt{\qone})\\
&~+\frac{1}{2n}\left(n\llbracket m\gam_1^2=2(m+1)\gam_2\rrbracket -1-\textcolor{red}{\eta(m\gam_1^2-2(m+1)\gam_2)}\sqrt{n}  \right)(-1)^mA_m(-n,-1/\sqrt{\qone}).
\end{align*}
\end{itemize}
\end{lemma}
\proof We first  note that only $y\in D$ contributes to $\Wb_m(\gam_1,\gam_2)$ because of our assumption $\gam_1,\gam_2\in D$. We may use the relation
\begin{align*}
\Wb_m(\gam_1,\gam_2)=\sum_{y\in D}\Vb_m(\gam_1-y,\gam_2-y^2)
\end{align*}
 and apply Lemma~1 to complete the proof. However, it is simpler to apply formula \eqref{eq:sieve} directly with
 \[
 h(x_1,\ldots,x_m)=\sum_{y\in D}\left\llbracket y+\sum_{j=1}^m x_j=\gam_1, y^2+\sum_{j=1}^m x_j^2=\gam_1^2-2\gam_2  \right\rrbracket.
 \]
 Since the argument is very similar to that used in the proof of Lemma~\ref{lemma1}, we just point out where the differences are.   Let $\va=(a_1,\ldots,a_t)$ be the vector of all cycle lengths of $\tau$ which are not multiples of $p$  (as defined  before), and let $\va'=(1,a_1,\ldots,a_t)$.
  Because of the extra variable $y$, there is no need to treat the case $l'(\tau)=0$ separately. We also note $b= m+1 \pmod p$ here, and the condition $b_0^2=ba_0$ becomes $m\gam_1^2=2(m+1)\gam_2$.
Therefore
\begin{align}
\frac{\Wb_m(\gam_1,\gam_2)}{m!}&=\frac{1}{m!}\sum_{\tau\in \Scal_m}(-1)^{m-l(\tau)} \qone^{l(\tau)-l'(\tau)}V_{1+l'(\tau)}(\va';\gam_1^2-2\gam_2,\gam_1)\nonumber\\
&=\frac{1}{m!}\sum_{\tau\in \Scal_m}(-1)^{m-l(\tau)} \qone^{l(\tau)-l'(\tau)}\left(\qone^{l'(\tau)-1}+R_{1+l'(\tau)}(\va';\gam_1^2-2\gam_2,\gam_1)\right)\nonumber\\
&=\frac{1}{\qone}{\qone\choose m} +(-1)^{m} \frac{1}{m!}\sum_{\tau\in \Scal_m}(-n)^{l(\tau)} \qone^{-l'(\tau)}R_{1+l'(\tau)}(\va';\gam_1^2-2\gam_2,\gam_1)  \label{eq:Wb8}\\
&=\frac{1}{\qone}{\qone\choose m}  \nonumber\\
&~+(-1)^{m}\left(\llbracket \gam_1=0,p\mid (m+1) \rrbracket \textcolor{red}{\eta(\gamma_2)}+\llbracket p\nmid (m+1)\rrbracket\frac{n\llbracket m\gam_1^2=2(m+1)\gam_2\rrbracket -1}{n}\right)\nonumber\\
&\quad \times\frac{1}{2}(A_m(-n,1/\sqrt{n})+A_m(-n,-1/\sqrt{n}))\nonumber\\
&~+(-1)^{m}\left(\llbracket \gam_1=0,p\mid (m+1) \rrbracket(n\llbracket \gam_2=0\rrbracket-1)+\llbracket p\nmid (m+1)\rrbracket \textcolor{red}{\eta( m\gam_1^2- 2(m+1)\gam_2)}\right)\nonumber\\
&\quad \times\frac{1}{2\sqrt{n}}(A_m(-n,1/\sqrt{n})-A_m(-n,-1/\sqrt{n})).  \label{eq:Wb9}
\end{align}
The proof is completed by  separating the cases $p\mid (m+1)$ and $p\nmid (m+1)$. \qed

\medskip

Substituting the values of $\Vb_{k+2}(\gam_1,\gam_2)$ and $\Wb_{k+1}(\gam_1,\gam_2)$ from Lemmas~1 and 2 into \eqref{eq:N2}, we immediately obtain the following generalization of \cite[Theorem~4.3]{ZhouWang17}.
\begin{thm}\label{thm:three}
Suppose $\qone$ is an even power of $p$, $D=\fqone$ is a subfield of $\fq$,  $\gam_1,\gam_2\in D$, and $A_m(u,w)$ be defined in \eqref{eq:Auw}. Then the following holds.
\begin{itemize}
\item[(a)] For $p\nmid (k+2)$,
\begin{align*}
&~~~~N(x^{k+2}+\gam_1x^{k+1}+\gam_2x^{k},r)\\
&=q^{k-r}{n\choose r}\sum_{j=0}^{k-r}(-q)^{-j}{n-r\choose r}+(-1)^{k-r}\frac{1}{n^2}{n\choose r}\left({n-r\choose k+2-r}-n{n-r\choose k+1-r}\right)\\
&~+\frac{(-1)^r}{2n}{k+1\choose r}\left(n\llbracket (k+1)\gam_1^2=2(k+2)\gam_2\rrbracket-1+\textcolor{red}{\eta((k+1)\gamma_1^2-2(k+2)\gamma_2)}\sqrt{n}\right)A_{k+1}(-n,1/\sqrt{\qone})\\
&~+\frac{(-1)^r}{2n}{k+1\choose r}\left(n\llbracket (k+1)\gam_1^2=2(k+2)\gam_2\rrbracket-1-\textcolor{red}{\eta((k+1)\gamma_1^2-2(k+2)\gamma_2)}\sqrt{n}\right)A_{k+1}(-n,1/\sqrt{\qone})\\
&~+\frac{(-1)^r}{2n\sqrt{n}}{k+2\choose r}\left(\textcolor{red}{\eta((k+1)\gamma_1^2-2(k+2)\gamma_2)}\sqrt{n}+n\llbracket (k+1)\gam_1^2=2(k+2)\gam_2\rrbracket-1\right)A_{k+2}(-n,1/\sqrt{\qone})\\
&~+\frac{(-1)^r}{2n\sqrt{n}}{k+2\choose r}\left(\textcolor{red}{\eta((k+1)\gamma_1^2-2(k+2)\gamma_2)}\sqrt{n}-n\llbracket (k+1)\gam_1^2=2(k+2)\gam_2\rrbracket+1\right)A_{k+2}(-n,-1/\sqrt{\qone}).
\end{align*}
\item[(b)] For $p\mid (k+2)$,
\begin{align*}
&~~~~N(x^{k+2}+\gam_1x^{k+1}+\gam_2x^{k},r)\\
&=q^{k-r}{n\choose r}\sum_{j=0}^{k-r}(-q)^{-j}{n-r\choose r}+(-1)^{k-r}\frac{1}{n^2}{n\choose r}\left({n-r\choose k+2-r}-n{n-r\choose k+1-r}\right)\\
&~-(-1)^{r+(k+2)/p}\frac{1}{n^2}{k+2\choose r}{n/p\choose (k+2)/p}\\
&~+\frac{(-1)^r}{2\sqrt{n}}{k+1\choose r}\llbracket \gam_1=0\rrbracket  \left(\textcolor{red}{\eta(\gamma_2)}\sqrt{n}+n\llbracket \gam_2=0\rrbracket-1\right)A_{k+1}(-n,1/\sqrt{\qone})\\
&~+\frac{(-1)^r}{2\sqrt{n}}{k+1\choose r}\llbracket \gam_1=0\rrbracket  \left(\textcolor{red}{\eta(\gamma_2)}\sqrt{n}-n\llbracket \gam_2=0\rrbracket+1\right)A_{k+1}(-n,-1/\sqrt{\qone})\\
&~+\frac{(-1)^r}{2n}{k+2\choose r}\llbracket \gam_1=0\rrbracket \left(n\llbracket \gam_2=0\rrbracket-1-\textcolor{red}{\eta(\gamma_2)}\sqrt{n}\right)A_{k+2}(-n,1/\sqrt{\qone}) \\
&~+\frac{(-1)^r}{2n}{k+2\choose r}\llbracket \gam_1=0\rrbracket \left(n\llbracket \gam_2=0\rrbracket-1+\textcolor{red}{\eta(\gamma_2)}\sqrt{n}\right)A_{k+2}(-n,-1/\sqrt{\qone}).
\end{align*}
\end{itemize}
\end{thm}

Substituting $(\gam_1,\gam_2)=(0,0),(1,0),(0,1),(1,1)$ into Theorem~\ref{thm:three}, we  immediately obtain the following corollaries.  Corollary~1 agrees with \cite[Theorem~4.3]{ZhouWang17} when $n=q$
by noting the following correspondence between the parameters:
\begin{align*}
k+2 &\longleftrightarrow n \\
r &\longleftrightarrow k\\
A_{k+2}(-q,1/\sqrt{q}) & \longleftrightarrow (-1)^n\alp(n)\\
A_{k+2}(-q,-1/\sqrt{q}) &\longleftrightarrow \beta(n).
\end{align*}

\begin{cor}\label{cor:one}
Suppose $\qone$ is an even power of $p$, $D=\fqone$ is a subfield of $\fq$.  Then the following holds.
\begin{itemize}
\item[(a)] For $p\nmid (k+2)$,
\begin{align*}
&~~~~N(x^{k+2},r)\\
&=q^{k-r}{n\choose r}\sum_{j=0}^{k-r}(-q)^{-j}{n-r\choose r}+(-1)^{k-r}\frac{1}{n^2}{n\choose r}\left({n-r\choose k+2-r}-n{n-r\choose k+1-r}\right)\\
&~+(-1)^r\frac{n-1}{2n}{k+1\choose r}\left(A_{k+1}(-n,1/\sqrt{\qone})
+A_{k+1}(-n,-1/\sqrt{\qone})\right)\\
&~+(-1)^r\frac{(n-1)\sqrt{n}}{2n^2}{k+2\choose r}\left(A_{k+2}(-n,1/\sqrt{\qone})-A_{k+2}(-n,-1/\sqrt{\qone})\right).
\end{align*}
\item[(b)] For $p\mid (k+2)$,
\begin{align*}
&~~~~N(x^{k+2},r)\\
&=q^{k-r}{n\choose r}\sum_{j=0}^{k-r}(-q)^{-j}{n-r\choose r}+(-1)^{k-r}\frac{1}{n^2}{n\choose r}\left({n-r\choose k+2-r}-n{n-r\choose k+1-r}\right)\\
&~-(-1)^{r+(k+2)/p}\frac{1}{n^2}{k+2\choose r}{n/p\choose (k+2)/p}\\
&~+(-1)^r\frac{(n-1)\sqrt{n}}{2n}{k+1\choose r}\left(A_{k+1}(-n,1/\sqrt{\qone})-
A_{k+1}(-n,-1/\sqrt{\qone})\right)\\
&~+(-1)^r\frac{n-1}{2n}{k+2\choose r}\left(A_{k+2}(-n,1/\sqrt{\qone})+A_{k+2}(-n,-1/\sqrt{\qone})\right).
\end{align*}
\end{itemize}
\end{cor}

\begin{cor}\label{cor:two}
Suppose $\qone$ is an even power of $p$, $D=\fqone$ is a subfield of $\fq$.  Then the following holds.
\begin{itemize}
\item[(a)] For $p\nmid (k+1)(k+2)$,
\begin{align*}
&~~~~N(x^{k+2}+x^{k+1},r)\\
&=q^{k-r}{n\choose r}\sum_{j=0}^{k-r}(-q)^{-j}{n-r\choose r}+(-1)^{k-r}\frac{1}{n^2}{n\choose r}\left({n-r\choose k+2-r}-n{n-r\choose k+1-r}\right)\\
&~+\frac{(-1)^r}{2n}{k+1\choose r}\left((\sqrt{n}-1)A_{k+1}(-n,1/\sqrt{\qone})-
(\sqrt{n}+1)A_{k+1}(-n,-1/\sqrt{\qone})\right)\\
&~+\frac{(-1)^r}{2n^2}{k+2\choose r}\left((n-\sqrt{n})A_{k+2}(-n,1/\sqrt{\qone})
+(n+\sqrt{n})A_{k+2}(-n,-1/\sqrt{\qone})\right).
\end{align*}
\item[(b)] For $p\mid (k+1)$,
\begin{align*}
&~~~~N(x^{k+2}+x^{k+1},r)\\
&=q^{k-r}{n\choose r}\sum_{j=0}^{k-r}(-q)^{-j}{n-r\choose r}+(-1)^{k-r}\frac{1}{n^2}{n\choose r}\left({n-r\choose k+2-r}-n{n-r\choose k+1-r}\right)\\
&~+\frac{(-1)^r(n-1)}{2n}{k+1\choose r}\left(A_{k+1}(-n,1/\sqrt{\qone})+A_{k+1}(-n,-1/\sqrt{\qone})\right)\\
&~+\frac{(-1)^r(n-1)\sqrt{n}}{2n^2}{k+2\choose r}\left(A_{k+2}(-n,1/\sqrt{\qone})
-A_{k+2}(-n,-1/\sqrt{\qone})\right).
\end{align*}

\item[(c)] For $p\mid (k+2)$,
\begin{align*}
&~~~~N(x^{k+2}+\gam_1x^{k+1}+\gam_2x^{k},r)\\
&=q^{k-r}{n\choose r}\sum_{j=0}^{k-r}(-q)^{-j}{n-r\choose r}+(-1)^{k-r}\frac{1}{n^2}{n\choose r}\left({n-r\choose k+2-r}-n{n-r\choose k+1-r}\right)\\
&~-(-1)^{r+(k+2)/p}\frac{1}{n^2}{k+2\choose r}{n/p\choose (k+2)/p}.
\end{align*}
\end{itemize}
\end{cor}

\begin{cor}\label{cor:three}
Suppose $\qone$ is an even power of $p$, $D=\fqone$ is a subfield of $\fq$.  Then the following holds.
\begin{itemize}
\item[(a)] For $p\nmid (k+2)$,
\begin{align*}
&~~~~N(x^{k+2}+x^{k},r)\\
&=q^{k-r}{n\choose r}\sum_{j=0}^{k-r}(-q)^{-j}{n-r\choose r}+(-1)^{k-r}\frac{1}{n^2}{n\choose r}\left({n-r\choose k+2-r}-n{n-r\choose k+1-r}\right)\\
&~+\frac{(-1)^r}{2n}{k+1\choose r}\left((\sqrt{n}-1)A_{k+1}(-n,1/\sqrt{\qone})-
(\sqrt{n}+1)A_{k+1}(-n,-1/\sqrt{\qone})\right)\\
&~+\frac{(-1)^r}{2n^2}{k+2\choose r}\left((n-\sqrt{n})A_{k+2}(-n,1/\sqrt{\qone})+
(n+\sqrt{n})A_{k+2}(-n,-1/\sqrt{\qone})\right).
\end{align*}
\item[(b)] For $p\mid (k+2)$,
\begin{align*}
&~~~~N(x^{k+2}+x^{k},r)\\
&=q^{k-r}{n\choose r}\sum_{j=0}^{k-r}(-q)^{-j}{n-r\choose r}+(-1)^{k-r}\frac{1}{n^2}{n\choose r}\left({n-r\choose k+2-r}-n{n-r\choose k+1-r}\right)\\
&~-(-1)^{r+(k+2)/p}\frac{1}{n^2}{k+2\choose r}{n/p\choose (k+2)/p},\\
&~+\frac{(-1)^r}{2n}{k+1\choose r}\left((n-\sqrt{n})A_{k+1}(-n,1/\sqrt{\qone})
+(n+\sqrt{n})A_{k+1}(-n,-1/\sqrt{\qone})\right)\\
&~-\frac{(-1)^r}{2n^2}{k+2\choose r}\left((\sqrt{n}+1)A_{k+2}(-n,1/\sqrt{\qone})
-(\sqrt{n}-1)A_{k+2}(-n,-1/\sqrt{\qone})\right).
\end{align*}
\end{itemize}
\end{cor}

\begin{cor}\label{cor:four}
Suppose $\qone$ is an even power of $p$, $D=\fqone$ is a subfield of $\fq$.  Then the following holds.
\begin{itemize}
\item[(a)] For $p\nmid (k+2)(k+3)$,
\begin{align*}
&~~~~N(x^{k+2}+x^{k+1}+x^{k},r)\\
&=q^{k-r}{n\choose r}\sum_{j=0}^{k-r}(-q)^{-j}{n-r\choose r}+(-1)^{k-r}\frac{1}{n^2}{n\choose r}\left({n-r\choose k+2-r}-n{n-r\choose k+1-r}\right)\\
&~+\frac{(-1)^r}{2n}{k+1\choose r}\left((\sqrt{n}-1)A_{k+1}(-n,1/\sqrt{\qone})-
(\sqrt{n}+1)A_{k+1}(-n,-1/\sqrt{\qone})\right)\\
&~+\frac{(-1)^r}{2n^2}{k+2\choose r}\left((n-\sqrt{n})A_{k+2}(-n,1/\sqrt{\qone})
+(n+\sqrt{n})A_{k+2}(-n,-1/\sqrt{\qone})\right).
\end{align*}
\item[(b)] For $p\mid (k+3)$,
\begin{align*}
&~~~~N(x^{k+2}+x^{k+1}+x^{k},r)\\
&=q^{k-r}{n\choose r}\sum_{j=0}^{k-r}(-q)^{-j}{n-r\choose r}+(-1)^{k-r}\frac{1}{n^2}{n\choose r}\left({n-r\choose k+2-r}-n{n-r\choose k+1-r}\right)\\
&~+\frac{(-1)^r(n-1)}{2n}{k+1\choose r}\left(A_{k+1}(-n,1/\sqrt{\qone})+
A_{k+1}(-n,-1/\sqrt{\qone})\right)\\
&~+\frac{(-1)^r(n-1)\sqrt{n}}{2n^2}{k+2\choose r}\left(A_{k+2}(-n,1/\sqrt{\qone})
-A_{k+2}(-n,-1/\sqrt{\qone})\right).
\end{align*}
\item[(c)] For $p\mid (k+2)$,
\begin{align*}
&~~~~N(x^{k+2}+x^{k+1}+x^{k},r)\\
&=q^{k-r}{n\choose r}\sum_{j=0}^{k-r}(-q)^{-j}{n-r\choose r}+(-1)^{k-r}\frac{1}{n^2}{n\choose r}\left({n-r\choose k+2-r}-n{n-r\choose k+1-r}\right)\\
&~-(-1)^{r+(k+2)/p}\frac{1}{n^2}{k+2\choose r}{n/p\choose (k+2)/p}.
\end{align*}
\end{itemize}
\end{cor}

\subsection{Asymptotic expressions for general $n$}

When $n$ is not an even power of $p$,   expressions for $\Vb_m(\gam_1,\gam_2)$ and  $\Wb_m(\gam_1,\gam_2)$  become complicated.  So we focus on the asymptotic estimate. For convenience, we shall use $\OT(A)$ to denote any number whose absolute value is less than or equal to $A$.

\begin{thm}\label{thm:four} Suppose  $D=\fqone$ is a subfield of $\fq$ and $\gam_1,\gam_2\in D$. Then
\begin{align*}
&~~N(x^{k+2}+\gam_1x^{k+1}+\gam_2x^k)\\
&=q^{k-r}{n\choose r}\sum_{j=0}^{k-r} (-q)^{-j}{n-r\choose j}\\
&~~+(-1)^{k+1-r}\frac{1}{\qone}{n\choose r}{\qone-r\choose k+1-r}+(-1)^{k+2-r}\frac{1}{\qone^2}{\qone\choose r}{\qone-r\choose k+2-r}\\\\
&~~+\OT\left({k+2\choose r}A_{k+2}(n,1/\sqrt{n})+\sqrt{n}{k+1\choose r}A_{k+1}(n,1/\sqrt{n})   \right).
\end{align*}
\end{thm}
\proof We apply the bound $\left|R_m(\va,a_0,b_0)\right|\le n^{-m/2}$ given in Proposition~\ref{prop:three}.
It follows from \eqref{eq:Stirling3} and \eqref{eq:Vbm} that
\begin{align*}
\frac{\Vb_{k+2}(\gam_1,\gam_2)}{(k+2)!}
&=\frac{1}{\qone^2}{\qone\choose k+2}+(-1)^{k}\frac{\qone^2\llbracket \gam_1=\gam_2=0\rrbracket -1}{\qone^2} A_{k+2}(-n,0)\\
&~~~+\OT\left(\frac{1}{(k+2)!}\sum_{\tau\in \Scal_{k+2},l'(\tau)>0}\qone^{l(\tau)-l'(\tau)/2} \right)\\
&= \frac{1}{\qone^2}{\qone\choose k+2}+(-1)^{k}\frac{\qone^2\llbracket \gam_1=\gam_2=0\rrbracket -1}{\qone^2}A_{k+2}(-n,0)\\
 &~~~+\OT\left(A_{k+2}(n,1/\sqrt{n})-A_{k+2}(n,0)\right).
\end{align*}
By \eqref{eq:Auw}, we have  $|A_{m}(u,v)|\le A_{m}(|u|,|v|)$, and hence
\begin{align*}
&~~~\left|(-1)^{k}\frac{\qone^2\llbracket \gam_1=\gam_2=0\rrbracket -1}{\qone^2} A_{k+2}(-n,0)+\left(A_{k+2}(n,1/\sqrt{n})-A_{k+2}(n,0)\right)\right|\\
&\le A_{k+2}(n,0)+\left(A_{k+2}(n,1/\sqrt{n})-A_{k+2}(n,0)\right)\\
&\le A_{k+2}(n,1/\sqrt{n}).
\end{align*}
It follows that 
\begin{align*}
\frac{\Vb_{k+2}(\gam_1,\gam_2)}{(k+2)!}
&=\frac{1}{\qone^2}{\qone\choose k+2}+\OT\left(A_{k+2}(n,1/\sqrt{n})\right).
\end{align*}

Similarly we use \eqref{eq:Wb8} to obtain
\begin{align*}
\frac{\Wb_{k+1}(\gam_1,\gam_2)}{(k+1)!}
&=\frac{1}{\qone}{\qone\choose k+1}+\OT\left(\frac{\sqrt{n}}{(k+1)!}\sum_{\tau\in \Scal_{k+1}}\qone^{l(\tau)-l'(\tau)/2} \right)\\
&=\frac{1}{\qone}{\qone\choose k+1}+\OT\left(\sqrt{n}A_{k+1}(n,1/\sqrt{n})   \right).
\end{align*}
Substituting the above two expressions into \eqref{eq:N2} and using \eqref{eq:BinIden}, we complete the proof. \qed

To illustrate that our error bound is exponentially smaller than the one given in \cite[Theorem~1.5]{LiWan20}, we consider the case $n=q=p^2$ and $k= p-3$.
Then we obtain from \eqref{eq:Auw} that
\begin{align*}
A_{k+2}(p^2,1/p)&={2p-2\choose p-1},\\
A_{k+1}(p^2,1/p)&={2p-3\choose p-2},
\end{align*}
and the error term given in Theorem~5 becomes
\begin{align*}
E(p,r)&={p-1\choose r}{2p-2\choose p-1}+p{p-2\choose r}{2p-3\choose p-2}\\
&=\left(1+\frac{p(p-1-r)}{2(p-1)}\right){p-1\choose r}{2p-2\choose p-1}.
\end{align*}
The corresponding error term given in \cite[Theorem~1.5]{LiWan20} becomes
\begin{align*}
E'(p,r)&={p-1\choose r}{4p-1\choose p-1}+p{p-2\choose r}{4p-2\choose p-2}\\
&=\left(1+\frac{p(p-1-r)}{4p-1}\right){p-1\choose r}{4p-1\choose p-1}.
\end{align*}
Noting $r\le k+2=p-1$ and using the well-known Stirling's formula
\[
p!\sim \sqrt{2\pi p}(p/e)^p,~~\hbox{ as } p\to \infty,
\]
we obtain
\[
\frac{E(p,r)}{E'(p,r)}\sim \sqrt{6}\frac{p+1-r}{p+3-r}(3/4)^{3p}, ~~\hbox{ as } p\to \infty.
\]

We also remark that our error bound is sharp when $k=o(\sqrt{p})$ because
\[
\lim_{p\to \infty} \frac{|A_m(-p^2,\pm 1/p)|}{A_m(p^2,1/p)}=1,~~\hbox{ when }m=o(\sqrt{p}).
\]

\section{The expected value and variance of the distance between a received word and a random codeword in $\RS_{n,k}$}
In this section,  let $D$ be any subset of $\fq$ with $|D|=n$.
The following bounds are known \cite{LiWan20}:
\begin{align*}
n-k-\ell\le d(f,\RS_{n,k})\le n-k.
\end{align*}

The upper bound above shows that there are codewords in $\RS_{n,k}$ which are within distance $n-k$ away from any given $f\in \Mcal_{k+\ell}$. The following result shows that the average distance between $f$
and codewords in $\RS_{n,k}$ is very close to $n$ and the variance is less than $n/q$.

\begin{thm}
Let $Z(f)$ denote the distance between a  received word represented by  $f\in \Mcal_{k+\ell}$
 and a random codeword in $\RS_{n,k}$ (under uniform distribution, that is, each word in $\RS_{n,k}$ is chosen with probability $q^{-k}$).
\begin{itemize}
\item[(a)] The expected value of $Z(f)$ is equal to
\begin{align*}
\frac{q-1}{q}n.
\end{align*}
\item[(b)]  The variance of $Z(f)$ is equal to
\begin{align*}
\frac{q-1}{q^2}n.
\end{align*}
\end{itemize}
\end{thm}
\proof We use the generating function approach as in Section~2.
We use another indeterminate $u$ to mark the number of occurrences of linear factors $x+\alp$ with $\alp\in D$. Let $G(z,u)$ denote the corresponding generating function. Then a similarly argument as in Section~2 gives
\begin{align*}
G(z,u)&=\prod_{\alp\in D}\left(\langle 1\rangle+u\sum_{j\ge 1}z^{j}\langle x+\alp\rangle^j\right)
\prod_{\beta\in \fq\setminus D}\left(\langle 1\rangle+\sum_{j\ge 1}z^j\langle x+\beta\rangle^j\right)
\sum_{g\in \Mcal'}z^{\deg (g)}\langle g\rangle \nonumber\\
&=\prod_{\alp\in D}\left(\frac{u}{\langle 1\rangle-z\langle x+\alp\rangle}+(1-u)\langle 1\rangle\right)\prod_{\alp\in D}(\langle 1\rangle-z\langle x+\alp\rangle)F(z)\nonumber\\
&=F(z)\prod_{\alp\in D}(\langle 1\rangle+(u-1)z\langle x+\alp\rangle),
\end{align*}
where $F(z)$ is given in \eqref{eq:F}, and $\Mcal'$ denote the subset of $\Mcal$ consisting of all  polynomials with no linear factors.

Thus
\begin{align*}
\frac{\partial}{\partial u}G(z;u)
&=zG(z,u)\sum_{\alp\in D}\frac{\langle x+\alp\rangle}{\langle 1\rangle+(u-1)z\langle x+\alp\rangle},\\
\frac{\partial}{\partial u}G(z;u)\Bigr|_{u=1}&=zF(z)\sum_{\alp\in D}\langle x+\alp\rangle,\nonumber\\
\frac{\partial^2}{(\partial u)^2}G(z;u)\Bigr|_{u=1}
&=z^2F(z)\left(\left(\sum_{\alp\in D}\langle x+\alp\rangle\right)^2-\sum_{\alp\in D}\langle x+\alp\rangle^2\right).
\end{align*}
Assuming $k\ge 1$ and using \eqref{eq:F}, we obtain
\begin{align*}
\left[\langle f\rangle z^{k+\ell}\right]\frac{\partial}{\partial u}G(z;u)\Bigr|_{u=1}
&=\left[\langle f\rangle z^{k-1}\right]\frac{q^{\ell}}{1-qz}E\sum_{\alp\in D}\langle x+\alp\rangle\nonumber\\
&=nq^{k-1},\\
\left[\langle f\rangle z^{k+\ell}\right]\frac{\partial^2}{(\partial u)^2}G(z;u)\Bigr|_{u=1}
&=\left[\langle f\rangle z^{k-2}\right]\frac{q^{\ell}}{1-qz}E\left(\left(\sum_{\alp\in D}\langle x+\alp\rangle\right)^2-\sum_{\alp\in D}\langle x+\alp\rangle^2\right)\nonumber\\
&=\left(n^2-n\right)q^{k-2}.
\end{align*}
Recall that $n-Z(f)$ is the number of linear factors of a polynomial counted by $G(z,u)$.

Since the total number of codewords in $\RS_{n,k}$ is equal to $q^k$, it follows that (see \cite[P.158]{FlaSed09})
the expected value and the variance of $n-Z(f)$ are, respectively,
\begin{align*}
\mu&=\frac{1}{q^k}\left[\langle f\rangle z^{k+\ell}\right]\frac{\partial}{\partial u}G(z;u)\Bigr|_{u=1}=\frac{n}{q},\\
\sigma^2&=\frac{1}{q^k}\left[\langle f\rangle z^{k+\ell}\right]\frac{\partial^2}{(\partial u)^2}G(z;u)\Bigr|_{u=1}+\mu-\mu^2=\frac{q-1}{q^2}n.
\end{align*}
Now the proof is completed by noting that the expected value of $Z(f)$ is equal to $n-\mu$, and
the variance of $Z(f)$ is equal to the variance of $n-Z(f)$. \qed

\section{Conclusion}

In this paper we used the generating function approach to find the number of codewords in $\RS_{n,k}$ which are distance $n-r$ away from
a given received word. This is equivalent to finding the number of monic polynomials with prescribed $\ell$ leading coefficients
and exactly $r$ distinct roots in a subset $D\subseteq \fq$.
The coefficients of the generating functions are from the group algebra $\mathbb{Q}\Ecal$ generated from the group
$\Ecal$ of equivalence classes of polynomials with prescribed leading coefficients. The coefficients of the generating functions
can be expressed in terms of the number of solutions to diophantine equations over $\fq$, where some variables are required
to have distinct values. The sieve formula of Li and Wan is applied to obtain simple formulas for the case of prescribing one leading coefficient ($\ell=1$)
and $D\cup \{0\}$ is an additive subgroup of $\fq$. This extends \cite[Theorem~3.1]{ZhouWang17}.
For the case of prescribing two leading coefficients ($\ell=2$), we are able to extend  \cite[Theorem~4.3]{ZhouWang17}
so that $D$ can be any subfield of $\fq$ and $|D|$ is an even power of an odd prime, and the two leading coefficients are arbitrary. When $|D|$ is a general power of an odd prime, we obtain
a simple asymptotic estimate which has much smaller error than the one given in \cite[Theorem~1.5]{LiWan20}. For the general subset $D$ (with no extra conditions),
we also find simple formulas for the first two moments of the distance between a random codeword in $\RS_{n,k}$ and a received word.
In a subsequent paper, we will apply our method to derive the corresponding asymptotic formulas for the case $\ell\ge 3$, which improve the error bound in \cite[Theorem~1.5]{LiWan20}.  The argument in Section~5 can also be used to derive simple expressions for higher moments.

\medskip

\centerline{Acknowledgement}

\medskip

I would like to thank Simon Kuttner for pointing out an error in Lemma~1 of the published version \cite{pub22} of this paper. This version fixes all the errors in section~4.1, which are due to the missing  quadratic character $\eta$ in Proposition~4.

\end{document}